\documentclass[10pt, oneside]{amsart}

\usepackage{courier}

\usepackage[T1]{fontenc}

\usepackage[%
	protrusion=true,
	expansion=false,
	auto=false,
	letterspace=100,
	]{microtype}

\usepackage{calc}
\usepackage{xcolor}
\usepackage{graphicx}
\graphicspath{{./figures/}}
\usepackage{tikz}
\usepackage[hang,flushmargin]{footmisc}
\usepackage{mathastext}
\usepackage[papersize={8.5in,11in},%
textwidth=5.5in,textheight=8.5in,%
headsep=2\baselineskip,ignorefoot,top=1.25in]{geometry}
\usepackage[normalem]{ulem}

\definecolor{linkred}{rgb}{0,0,0}
\definecolor{linkblue}{rgb}{0,0,0}
\PassOptionsToPackage{hyphens}{url} 
\usepackage[
    setpagesize=false,
	pdfpagelabels=false,
    pdfstartview={Fit 1000},
    bookmarksnumbered=false,
    linktoc=all,
    colorlinks=true,
    anchorcolor=black,
    menucolor=black,
    runcolor=black,
    filecolor=black,
    linkcolor=linkblue,
	citecolor=linkblue,
	urlcolor=linkred,
]{hyperref}

\newcommand\mfg{\text{\large$\displaystyle{\mathfrak{g}}$}}
\newcommand{\GG}{\text{\small$\displaystyle{\mathbb{G}}$}}
\newcommand{\QQ}{\text{\small$\displaystyle{\mathbb{Q}}$}}
\newcommand{\ooplus}{\text{\Large$\displaystyle{\,\oplus\,}$}}
\newcommand{\ccap}{\text{\Large$\displaystyle{\,\cap\,}$}}
\newcommand{\ccup}{\text{\Large$\displaystyle{\,\cup\,}$}}
\newcommand{\nnum}[1]{\refstepcounter{equation}\label{#1}\theequation}
\newcommand{\kj}{\vskip0.8\baselineskip}
\newcommand{\kjm}{\vskip0.4\baselineskip}
\newcommand\blfootnote[1]{%
  \begingroup
  \renewcommand\thefootnote{}\footnote{#1}%
  \addtocounter{footnote}{-1}%
  \endgroup
}
\def\limit{\qopname\relax m{limit~}}
\def\qed{{\null\hfill{}Q.E.D.}}
\def\Hom{\mathrm{Hom}}

\raggedright
\makeatletter
\DeclareMathSizes{\f@size}{\f@size}{\f@size}{\f@size}
\def\<{\kern4pt\(}
\def\>{\)%
\@ifnextchar{-}{}{
\@ifnextchar{.}{}{
\@ifnextchar{,}{}{
\@ifnextchar{;}{}{
\@ifnextchar{:}{}{
\@ifnextchar{'}{}{
\@ifnextchar{)}{}{
\kern4pt{}
}}}}}}}
}
\def\ps@headings{\ps@empty
\def\@evenhead{%
\setTrue{runhead}
\normalfont
\rlap{\thepage}\hfil
\leftmark{}{}\hfil}\bigskip%
\def\@oddhead{%
\setTrue{runhead}%
\normalfont \hfil
\rightmark{}{}\hfil \llap{\thepage}}%
\let\@mkboth\markboth
}
\makeatother
\newcommand{\kproc}[1]{\noindent\uline{#1}}
\begin{document}

\centerline{\large Instability in invariant theory}
\centerline{by}
\centerline{\large George Kempf}\blfootnote{\raggedright\normalsize{}\vskip-18pt\noindent{}Research partially supported by NSF Grant MPS75-0S578

\noindent\linespread{1}\selectfont{}Transcriber's note: In the fall of 1976, my advisor, David Mumford, handed me a short preprint by George Kempf to read. It was the first state of what eventually became his influential Annals paper ``Instability in Invariant Theory'' (Annals of Mathematics, Second Series, Vol. 108, No. 2 (Sep., 1978), pp. 299-316). The introduction to the published version ends with an acknowledgement and a dig: ``I~want to thank the referee of the first version of this paper for pointing out Corollary 4-5 and conjecturing that the original \< \{0\} \>-instability could be replaced by \< S \>-instability. Unfortunately, the inclusion of \< S \>-instability has completely destroyed the simplicity of the original version.'' Over the intervening years, the simplicity and elegance of the first version has continued to create a readership for it, and copies (see Figure~\ref{8}) continue to circulate informally. I created this LaTeX'ed version to make the paper accessible to all who may be interested in it, trying to keep the look close to that of the original typewritten preprint and making changes only to correct a few obvious typos and harmonize the markup. My thanks to George's children, Robin and Lucas Kempf, for graciously granting me permission to post this transcription.\\[3pt]Comments welcome to Ian Morrison (morrison@fordham.edu).}

\thispagestyle{empty}
\enlargethispage{\baselineskip}
Let \< v \> be a point of the representation space \< V \> of a reductive algebraic group \< G \>. The point \< v \> is called unstable if any polynomial function 
on \< V \>, which is invariant under \< G \> and vanishes at \< 0 \> in \< V \>, must vanish at 
\< v \>. There is a related notion of instability. Let \< \lambda:\GG_m \rightarrow G \> be a one-parameter subgroup of G. The point \< v \> is called \< \lambda \>-unstable if 
\<\displaystyle{\limit_{\kern-2pt{}t\text{\scriptsize$\to$} 0}} \lambda(t)\cdot v = 0 \>. 

The Hilbert-Mumford numerical criterion for unstability states that \< v \> 
is unstable if and only if \< v \> is \< \lambda \>-unstable for some one parameter subgroup 
\< \lambda \> of \< G \>. This criterion was given in Chapter 2 of~\cite{4} for linearly reductive 
groups. Subsequently, it has been established for arbitrary reductive groups 
\cite{5,6} due to the efforts of C.~S.~Seshadri, M.~Nagata and W.~Haboush among others. 
 
Given an unstable point \< v \>, Mumford has posed the problem of picking 
out a natural class \< \Lambda_v \> of one-parameter subgroups \< \lambda \> of \< G \> such that \< v \> 
is \< \lambda \>-unstable. I will present a solution to this problem. Surprisingly, this 
result requires very little not contained in Mumford's book. 

I don't know any solution to Tit's more general ``center'' conjecture, 
but my result is strong enough to establish the generalization of Godement's 
conjecture mentioned by Mumford on page 64 of his book. For the case of the 
real ground field, this result has been proven by D. Birkes~\cite{1}. 

There remains the question of providing some intuitive geometric 
characterization of the class \< \Lambda_v \> and investigating the properties of the ``flag of highest contact'' \< P_v \> (see page 48 of~\cite{4}). For this reason, I have 
translated much of Mumford's treatment into concrete representation theory 
terminology. Another form of the Tits building was presented by Mumford in~\cite{3}.

\kj

\kproc{\S}1. Let \< M \> be a finite dimensional real vector space with positive definite 
inner product \<  (~,~ ) \>. Let \< ||m|| = (m,m)^{1/2}\> for \< m \in  M \>. Let \< S \> be the 
sphere \<  ||m|| = 1 \> in \< M \>. 

Let \< F \> be a finite set of real-valued linear functions on \< M \>. For
each point \< m \> of \< M \>, set \< h(m) \> equal to the minimum value of \< f(m) \> for all 
\< f  \> in \< F \>. We now state the
 
 \kj
 
\kproc{Lemma} \nnum{1}. Assume that the function \< h \> actually has a positive value somewhere 
in \< M \>. Let \< U = \{m \mid h(m) > O\} \>. There is a unique point \< p_h \> of \< S \ccap U \>, where 
\< h \> obtains a maximum value. In fact, \< p_h \> is the only point of \< S \ccap U \>, where 
\< h \> has a relative (local) maximum value. 

\kj

\kproc{Proof}. Clearly, \< h \> is a continuous function and must obtain a maximal value 
when it is restricted to the compact sphere \< S \>. Let \< p \> be any point of \< S \>, where 
\< h \> obtains this maximal value. By assumption, \< p \in S \ccap U := U'  \>and \< h(p) > O \>. 

Assume that \< M \> is a line. Then, \< S = \{p\} \ccup \{-p\} \>. By definition, 
\< h(p) \le -h(-p) \>. Thus, \< h(p) > 0 > h(-p) \> and \< U'= {p} \>. This case is trivial. 

A more interesting case is when \< M \> is a plane. The sets \< U_f = \{s \in Sj \mid f(s) > O\} \> 
for \< f \in F \> are open half-circles and \< U':=\ccap U_f \> is an 
open arc. To prove the lemma in this case, we note that \< f \> restricted to 
\< U_f  \> is a strictly convex function (with respect to angle ordering on \< U_f  \>). 
Thus, \< h \> is a strictly convex function on \< U \>. and must therefore have its 
sole relative maximum value at \< p \>. 
The general case follows by remarking that, for any other point \< s \> of 
\< S \>, \< s \> and \< p \> either span a line or a plane. \qed
 
 \kjm
 
 Let \< R_h \> denote the open ray through \< p_h \>: \( R_h \>  is the ray along which 
 \< h \> increases most rapidly. 
 
 We will also need an integrality statement about this ray. Let \< L \> 
 be a lattice in \<V \> . We will assume that the inner product of two elements 
 of \< L \> is integral. Furthermore, we assume that each function in \< F \>  has 
 integral values on \< L \> . 

 With all the above assumptions, we have 
 
 \kj

\kproc{Lemma} \nnum{2}. The ray \< R_h\> contains an element \< \lambda_h \in L\>  such that any element in \< R_h \ccap L\> is a positive integral multiple of \< \lambda_h \in L\>.

\kj 

 \kproc{Proof}. If we show that \< R_h \ccap L\> is not empty, the statement follows because the 
 intersection of \<L\> with the line generated by  \< R_h\>  will be a rank one abelian 
 group. Let \<G\> be the subset of \<F \> consisting of the functions \<f\> such that 
 \<f(p_h) = h(p_h)\>. \(M' = \{m \in M \mid g(m) = g'(m)\>~ for all \<g\> and \<g\>' in \,\<G\}\>. Let 
\< F'\> and \<h'\> be the restrictions of \< F\> and \<h\> to \<M'\>. Clearly, \<p_h = p_{h'}\>
 as \<p_h\> is a positive maximum for \<h'\>. Furthermore, \<L' = M' cap L\> n is a lattice 
 in \<V'\> as the equations of \<M'\> are integral with respect to \<L\>. The lemma 
 for \<h'\> clearly implies the lemma for \<h\>. 

 Therefore, we may assume that there is a function \<g\> in \<F\> such that 
\<g(p+h)=h(p_h) > f(p_h)\>  for any other function \<f\> in \<F\>. Thus, \<p+h\> must be 
a positive relative maximum for \<g\> as \<g = h\> near to \<p+_h\>. Hence, it will be 
 enough to prove the lemma when \<F = \{g\}\>. 

 We may end the proof by noting that, if \<g^*\> is the point of \<M\> such 
 that \<(g^*, m) = g(v)\>, then, \<g^* \in \QQ\cdot L\> by the integrality of \<( ~, ~)\>. Furthermore, \< g^*/||g^*||\> is evidently equal to \<p_h\> in this case. \qed

\kjm

If one dropped the assumption that \<h\> takes a positive value, the 
above argument shows 
\kjm
\kproc{Lemma} \nnum{3}. The function \<h\> has at least one point on the sphere, where it takes 
a maximal value. Some ray in the direction of a point on the sphere where \<h\> 
has a relative maximum must contain an element of \<L\> when we have the integrality assumption. 

\kj

\kproc{\S2}. Let \< V \> be a finite dimensional representation of a reductive group \< G \>. 
Let \< S \> be any torus of \< G \>. We have a weight decompositon, \< V =\ooplus V^\chi \>. The 
\< V^\chi\> are non-zero eigenspaces for distinct characters \< \chi \> of \< S \>. These 
characters are called the \< S \>-weights of V and the \< V^\chi \>'s wi11 be called 
weight spaces. 

Let \< v \> be any non-zero element of \< V \>. The state of \< v \> with respect 
to \< s \> is the set of weights \< \chi \> such that the projection of \< v \> onto \< V^\chi \>
is non-zero. Clearly, the state of a vector may be any arbitrary non-empty 
subset of the weights of \< V \>. 

Let \< \lambda: \GG_m\rightarrow G \> be a one-parameter subgroup of \< G \>. Let \< \chi \> denote a character of the image of \< \lambda \>. We have \<\chi(\lambda(t))\equiv t^{\chi(\lambda)}\>
for all \< t \> in \< \GG_m \> where \<  \chi(\lambda) \> is the integer defined by the formula. 
 Thus, the characters of \<Im(\lambda)\> and, hence, its eigenspaces in \< V \> are linearly ordered by the integer \< \chi(\lambda) \>. Define \< m(v,\lambda) \> to be the minimum \< \chi(\lambda) \>, where \< \chi \> runs through the state of \< v \> with respect to \<Im(\lambda)\>.

This last integer may be used to determine when \< v \> is \< \lambda \>-unstable, i.e.
\<\displaystyle{\limit_{\kern-2pt{}t\text{\scriptsize$\to$} 0}} \lambda(t)\cdot v = 0 \>. 
In fact, it is clear that \< v \> is \< \lambda \>-unstable if and only 
the integer \<  m(v,\lambda) > 0 \>. 

I want to record a simple property of the integer \<  m(v,\lambda) \>
\kjm
\noindent(*) \< \kern72pt m(v,\lambda) = m(g\cdot v, g \lambda g^{-1})\> for any element \<g\> of \< G\>.\\
\kjm
\noindent{}In fact, this is consequence of the stronger statement that the numerical 
invariants \<\chi(\lambda\}\> of the state of \< v \> with respect to \<Im(\lambda)\> are the same 
if we replace \<\lambda\> by \< g \lambda g^{-1}\> and \<v\> by \< g \cdot v\>.

Recall that the one-parameter subgroups of a torus \<  T \> are a free 
abelian group of rank equal the dimension of \< T \>. One defines a notion of 
length \< ||\lambda|| \> to any one-parameter subgroup \<\lambda\> of \<  G \> such that 
\kjm
(a)  \( ||g \lambda g^{-1}|| =||\lambda|| \>  for all \< g \> in \< G \> and any one parameter 
\\\leavevmode\kern21pt{}subgroup \<\lambda\> of \< G \>, and
\kjm
(b) for any maximal torus \< T \> of \< G \>, there is a bilinear positive 
\\\leavevmode\kern21pt{}definite integral-valued bilinear form \< (~,~) \> on the group of 
\\\leavevmode\kern21pt{}one-parameter subgroups\<\lambda\>. of \< T \> such that \< (\lambda,\lambda)^{1/2}= ||\lambda|| \>.
\kjm

\noindent{}Note that (a) implies that the inner product  \< (~,~) \> must be invariant under the 
Weyl group of \< G \> with respect to \< T \>. Conversely, given any invariant pairing, 
it extends to a unique notion of length of arbitrary one-parameter subgroups 
of \< G \>. 

With a fixed notion of length we will give a form of Proposition 2.17 
of Mumford's book. The result is 

\kj
\kproc{Lemma} \nnum{4}. Let \< v \> be a non-zero element of \< V \>. Let 
\kjm
\noindent\kern90pt \(\displaystyle{B(v) \equiv \sup \frac{m(v, \lambda)}{||\lambda||}}\)
\kjm
\noindent{}for all one-parameter subgroups \< \lambda \> of \< G \>. Then there exists at least one 
\< \lambda_0\> such that 

\noindent\kern90pt \(\displaystyle{B(v) =\frac{m(v, \lambda_0)}{||\lambda_0||}}\)~.
\kj
\kproc{Proof}. By the conjugacy properties (*) and a), we need only prove that, there 
is a one-parameter subgroup  \< lambda_0  \> of a fixed maximal torus \< T \> such that 

\noindent\kern90pt \(\displaystyle{\frac{m(g\cdot v, \lambda_0)}{||\lambda_0||} \ge \frac{m(g\cdot v, \lambda)}{||\lambda||}}\)
\kjm
\noindent{}for all \< g \> in \< G \> and \< \lambda \> in \< \Hom(\GG_m,T) \equiv \Gamma(T) \>, as any one-parameter subgroup of \< G \> is conjugate to an element of \<  \Gamma(T) \>. 

Let \< R \> be the state of a vector \< v' \> in \< V \> with respect to \< T \>. For 
any \< \lambda \>. in  \<  \Gamma(T) \>, we have \< m(v',\lambda) = \min \chi(\lambda) \> for \< \chi \> in the state \< R \>. 
As the \< \chi(\lambda)\> are integral-valued linear functions of \< \lambda \> in \<  \Gamma(T) \>, we may apply Lemma~\ref{3} to find a one-parameter subgroup \< \lambda_R \> in  \<  \Gamma(T) \> such that 
\kjm
\noindent\kern90pt \(\displaystyle{\frac{m( v', \lambda_R)}{||\lambda_R||} \ge \frac{m( v', \lambda)}{||\lambda||}}\)
\kjm
\noindent{}for any \< \lambda \> in \<  \Gamma(T) \>. 

The lemma now follows from the remark that these are only a finite 
number of possible states for vectors in \< V \).\qed
\kj

\kproc{\S3}. We next recall the parabolic subgroup \< P (\lambda) \> of \< G \> associated to a one-parameter subgroup \< \lambda \> of \< G \>. Consider the adjoint representation of \< G \> on 
its tangent space \<\mfg\>  at the identity. The subgroup \< P (\lambda) \> may be characterized 
as the subgroup of \< G \>, whose tangent space at the identity consists of the 
elements \< D \> of \<\mfg\>  such that \< m(D ,\lambda) \ge 0 \> plus the zero vector of  \<\mfg\>.

\( P (\lambda) \> has a marked Levi-subgroup \< L (\lambda) \> where \< L (\lambda) \> is the connected 
component of the elements of \< G \> which commute with \< \lambda \>. The tangent space of 
the unipotent radical \< U(\lambda) \> is the space of \< \lambda \>-unstable elements of \<\mfg\>.

Let \< V =\ooplus V^\chi \> be the \< \Im(\lambda) \>-weight decomposition of \< V \>. As \< \Im(\lambda) \> 
is contained in the center of  \< \lambda \> the action of  \< \lambda \> on \< V \> must be pre-
served by the weight spaces \< v^\chi \>.  In general, the action of \< P (\lambda) \> only preserves 
the weight filtration of \< V \>. 

Let \< V^i =\ooplus V^\chi \> 
for all weights \<\chi\> such that \<\chi(\lambda) \ge\> some integer i. 
Then, \< V^{i+1}\subset   V^i \> and the \< V^i\>'s  form a filtration of \< V \>. From the definition 
of \< P (\lambda) \>, one may check that the action of \< P (\lambda) \>  on \< V \> must preserve the weight  filtration \< V^i\>. In fact, the unipotent radical \< U(\lambda) \> acts trivially on the 
quotients \< V^i/ V^{i+1}\>. 

The next lemma is related to Mumford's Proposition 2.7. 

\kj 
\kproc{Lemma} \nnum{5}. Let \< v \> be a non-zero vector in \< V \>. Then, 

(a) \(m(v, \lambda) = \max i\> such that \<v \in  V^i\>. 

(b)  \(m(v, \lambda) =  m(p\cdot v, \lambda)\>  for any \< p \> in  \< P (\lambda) \>.

(c)  \(m(v, \lambda)\le  m(v', \lambda)\> for any non-zero vector \< v' \> in the closure of
\\\leavevmode\kern21pt{}the \< P (\lambda) \>-orbit of \<  v \>. 

\kj

\kproc{Proof}.\kern2pt{}(a) follows directly from the definition of the  \< V^i\>'s.

(b) is implied by (a) and the \< P (\lambda) \>-invariance of the \< V^i\>'s.

(c) follows because \< V^i\> contains the closure of \< V^i\setminus V^{i+1}\>. \qed

\kjm

We now are in a position to understand the central result of this paper. The statement will use the language introduced in and for Lemma~\ref{4}. 

\kj

\kproc{Theorem} \nnum{6}. Let \< v \> be a non-zero unstable vector in \< V \>. Let \< \Lambda_v \> be the set of  all one-parameter subgroups \< \lambda \> of \< G \> such that \<\displaystyle{\frac{m(g\cdot v, \lambda)}{||\lambda||} = B(v)}\> and  \< \lambda \> is not divisible as a subgroup of \< G \>. Then, 

(a)  \( \Lambda_v \> is not empty.

(b) there is a parabolic subgroup \< P_v \> of \< G \> such that \< P_v = P (\lambda) \>  for
\\\leavevmode\kern21pt{}any \< \lambda \> in  \< \Lambda_v \>, and 

(c)  \( \Lambda_v \> is a full conjugacy class of one-parameter subgroups of \< P_v \>. 

\kj

\kproc{Proof}. By Lemma~\ref{4}, there is at least one one-parameter subgroup \< \lambda_0 \> of \< G \>, where \<\displaystyle{\frac{m(g\cdot v, \lambda)}{||\lambda||}}\> obtains its maximum \< B(v)\>.
We may assume that  \< \lambda_0 \> is not 
divisible as a subgroup of \< G \>. Thus,\< \lambda_0 \in  Lambda_v\> and the statement (a) is true.
In this case, \< B(v) \> is positive by the Hilbert-Mumford criterion as \< v \> is 
unstable. 

Let \<T_0\> be any maximal torus containing \< \lambda_0 \>. If we use the reasoning
of Lemma~\ref{4}, then Lemmas~\ref{1} and~\ref{2} show that  \< \lambda_0 \> must be the only one-parameter subgroup of \<T_0\> where \<\displaystyle{\frac{m(g\cdot v, \lambda)}{||\lambda||} = (v)}\>. Furthermore, any maximal torus \< T \> of  \< P (\lambda_0) \>
is conjugate to \< T_0 \> by an element \< p \> of \< P (\lambda_0) \>; i.e. \< T = p^{-1}T_0 \,p) \>.
By Lemma~\ref{5} and (*), \<m(v, \lambda_0) = m(p\cdot v, \lambda_0)=m(v, p^{-1}\lambda_0 \,p)\>.
Thus, we may conclude that \< T \> contains a unique element \< p^{-1}\lambda_0 \,p \> of \< \Lambda_v \>. Clearly, \< P(\lambda_0)=  p^{-1} P(\lambda_0)p = P(p^{-1}\lambda_0 \,p) \>.

To finish the proof, let \<\lambda_1\> be another member of \< \Lambda_v \>. The intersection 
\< P(\lambda_0)\ccap  P(\lambda_1) \> of these two parabolic subgroups must contain a maximal 
torus \< T \> of \< G \>. Let  \<\lambda\> be the unique subgroup of \< T \> in  \< \Lambda_v \>. Then,
\< P(\lambda_0) = P(\lambda) = P(\lambda_1)\> by the last paragraph. This is statement (b). The statement  c) follows because we have seen that  \<\lambda_0\>,  \<\lambda\> and  \<\lambda_1\> are all conjugate in   \< P_v \>. \qed

\kj

\kproc{\S4}. In this section, we will give a rationality consequence of the Theorem~\ref{6}. 
Fix a perfect field \( k \). We will assume that \( G \) and its representation on \( V  \)
are all defined over \( k \). Furthermore, \( v \) will be a non-zero \( k \)-rational vector 
in \( V \). 

The statement of the next result requires that I explain what it means 
for the length function \( ||~~|| \) to be defined over \( k \). Let \( \sigma \) be an automorphism 
of the algebraic closure \( k' \) of \( k \) which fixes \( k \). Then, a length function 
\( ||~~|| \) is defined over \( k \) if

\kjm

\noindent(\#)\kern24pt{}\( ||^{\kern1pt{}\sigma\kern-2pt}\lambda|| =  ||\lambda||\) for all one-parameter subgroups \<\lambda\>  of \( G \) and all \<\sigma\>.

\kjm

To convince the reader of the existence of lengths defined over \( k \), 
we may take a maximal torus \( T \) of \( G \), which is defined over \( k \)~\cite[Theorem~18.2]{2}.
Any length function defined over \( k \) comes from an inner-product on the one-parameter 
subgroups \<\Gamma(T)\> of \( T \) which is integral-valued and positive-definite, 
which is invariant under the Weyl group of \( G \) over \( T \) and satisfies the 
equation (\#) for any \<\lambda\> in \( T \). As \( T \) is split over a finite extension of \( k \), 
the galois group of \( k'/k \) acts on \<\Gamma(T)\> by a finite quotient. Thus, there is never 
any problem satisfying (\#) if we sum any inner product over this finite group. 
The result of this section is 

\kj

\kproc{Theorem} \nnum{7}. Assume further that our length function is defined over k. Then, 

(a)  \( \Lambda_v \> is invariant under the galois group of \( k'/k \).

(b) \( P_v \> is defined over \( k \),

(c) there is a one-parameter subgroup \<\lambda\> in \< \Lambda_v \> which is defined 
\\\leavevmode\kern21pt{}over \( k \). 

\kj

\kproc{Proof}. Let \<\sigma\> be an element of the galois group of \( k'/k \). Then, 
\< m(v,\lambda) = m(v,^{\kern1pt{}\sigma\kern-2pt}\lambda)\> for any one-parameter subgroup \<\lambda\>  of \( G \) as \( v \) and the action of \( G \) are defined over \( k \). Thus, 
\<\displaystyle{\frac{m(g\cdot v, \lambda)}{||\lambda||}= \frac{m(g\cdot v, ^{\kern1pt{}\sigma\kern-2pt}\lambda)}{||^{\kern1pt{}\sigma\kern-2pt}\lambda||}}\> because  \( ||~~|| \) is defined 
over \( k \). Hence, the statement (a) follows from the definition of \< \Lambda_v \>. 

By the statement (b) of Theorem~\ref{6}, we have \< P(\lambda)=  P_v= P(^{\kern1pt{}\sigma\kern-2pt}\lambda_0)\> for any 
any element \<\lambda\> of \< \Lambda_v \>.  By the definition of the subgroup \<P(\lambda) \>, we have \<P(^{\kern1pt{}\sigma\kern-2pt}\lambda) = ^{\kern1pt{}\sigma\kern-2pt}P(\lambda) \>. Thus, \<P_v = ^{\kern1pt{}\sigma\kern-2pt}P_v \> for any \<\sigma\>. This proves (b).

For c), take a maximal torus \( T \) of \<P_v\> which is defined over \( k \). 
Then, \( T \) is also a maximal torus of \(  G \) as \<P_v\> is a parabolic subgroup. In the 
proof of Theorem~\ref{6}, we have seen that \( T \) has a unique one-parameter subgroup 
contained in  \< \Lambda_v \>. By uniqueness, this subgroup must be fixed by any  \<\sigma\> and,
hence, it is defined over \( k \). This proves c). \qed

\kj\kj\kj

\renewcommand{\refname}{}
\centerline{\uline{References}}\vskip-\baselineskip\vskip-\baselineskip

\bibstyle{plain}

\vfill\eject

\vspace*{-120pt}\leavevmode\kern-144pt\includegraphics[scale=1]{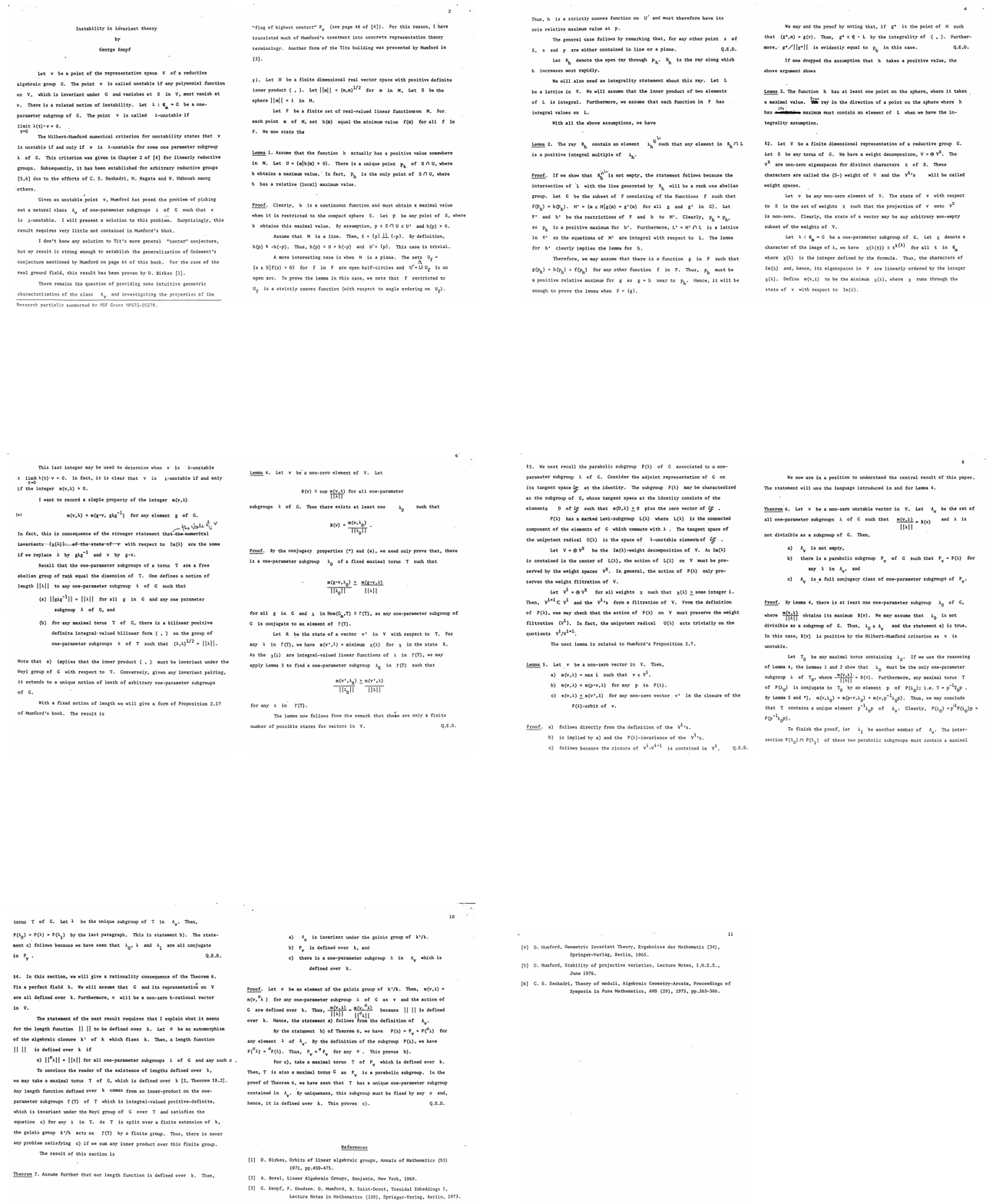}

\vspace*{-48pt}
\centerline{\kproc{Figure} \nnum{8}. A poorly scanned version of Kempf's original preprint.}

\end{document}